\documentclass[12pt]{amsart}

\usepackage{amsmath,amsthm,amssymb, pgfplots, tikz, color}
\usetikzlibrary{patterns}

\usepackage[margin=1in]{geometry}
\usetikzlibrary{positioning}

\newcommand{\Z}{\mathbb{Z}}

\newcommand{\Q}{\mathbb{Q}}

\newcommand\Disc{\operatorname{Disc}}

\newcommand\GL{\operatorname{GL}}
\newcommand\Sym{\operatorname{Sym}}

\numberwithin{equation}{section}

\newtheorem{theorem}{Theorem}

\theoremstyle{remark} \newtheorem*{remark}{Remark}

\numberwithin{theorem}{section}

\author{Robert J. Lemke Oliver}
\thanks{The first author was supported by an NSF Mathematical Sciences Postdoctoral Research Fellowship at Stanford University}
\address{Department of Mathematics, Tufts University, 503 Boston Ave, Medford, MA 02155}
\email{robert.lemke\_oliver@tufts.edu}
\author{Frank Thorne}
\thanks{The second author is partially supported by the National Science Foundation under Grant No. DMS-1201330.}
\address{Department of Mathematics, University of South Carolina, 1523 Greene St, Columbia, SC 29201}
\email{thorne@math.sc.edu}

\title{The number of ramified primes in number fields of small degree}

\begin{document}

\begin{abstract}
In this paper we investigate the distribution of the number of primes which ramify in number fields of degree $d \leq 5$. In analogy with the classical 
Erd\H{o}s-Kac theorem, we prove for $S_d$-extensions that the number of such primes is normally distributed with mean and variance $\log\log X$.  
\end{abstract}

\maketitle

\section{Introduction and statement of results}

Consider the following problem: How many number fields are there of absolute discriminant less than $X$?  By work of Minkowski and Hermite, 
this number is known to be finite, and so it is natural question to ask for an asymptotic estimate. 

Hermite's methods suggest that it is natural to count number fields of each fixed degree separately. To this end, write
\begin{equation}\label{eq:def_fd}
\mathcal{F}_d(X) := \{K/\mathbb{Q} : \mathrm{deg}\, K = d, |D_K|\leq X\}, \ \ \ \ N_d(X) := \# \mathcal{F}_d(X),
\end{equation}
where $D_K$ is the discriminant of $K$. It has been proved for each $d \leq 5$ that $N_d(X) \sim c_d X$ for an appropriate
constant $c_d$; the results are classical for $d = 2$ and due to
Davenport and Heilbronn \cite{DH1971} ($d = 3$), Cohen, Diaz y Diaz, and Olivier \cite{CDO2002, Cohen2003} and 
Bhargava \cite{Bhargava2005} ($d=4$), and Bhargava \cite{Bhargava2010} ($d = 5$). For higher $d$ this has been conjectured by 
Malle \cite{Malle2004} (see also Bhargava \cite{Bhargava2007}), and an upper bound roughly of the form $X^{d^\epsilon}$ was proved by Ellenberg and Venkatesh \cite{EV2006}.

In the cases where we have an asymptotic formula, one may also ask about the rate of convergence for these formulas. Write 
$\mathcal{F}_{S_d}(X)$ and $N_{S_d}(X)$ for the analogues of \eqref{eq:def_fd} where only those fields whose Galois closure is
an $S_d$-extension of $\mathbb{Q}$ are enumerated. Then, for $d \leq 5$, power saving error terms
have been proved for $N_{S_d}(X)$, by Belabas, Bhargava, and Pomerance \cite{BBP2010} ($d = 3, 4$)
and Shankar and Tsimerman \cite{ST2014} $(d = 5)$. 

It is also possible to count number fields satisfying a given finite set of local conditions, e.g., number fields in which a fixed prime $p$ is ramified.
The error terms will depend on the local conditions, but this dependence can be explicitly controlled, and it is again possible to give power saving error terms. 
Such results are essentially known in complete generality, although they are not all stated explicitly enough in the literature.\footnote{As we note later, such results
did appear after we completed the initial version of our paper. The strongest results to date are due to Ellenberg, Pierce, and Wood \cite{EPW}.}
In Section \ref{sec:field-counting} we give
explicit results of this nature which suffice to prove our main theorem.

\smallskip
It is typical that good estimates for the distribution of an arithmetic sequence in arithmetic progressions open the door to
further distribution results. For example:
\begin{itemize}
\item
Yang \cite{Yang2009}, and independently Cho and Kim \cite{CK2014}, used results on cubic fields to obtain the distribution of the low-lying
zeroes of the associated Dedekind zeta functions. Yang also did the same for quartic fields, and Cho and Kim did the same for some particular 
families of higher degree number fields (conditionally on several unproved hypotheses).
\item
In a second paper \cite{CK2014_CLT}, Cho and Kim prove a central limit theorem
for the partial sums of coefficients of Artin $L$-functions associated to $S_d$-number fields with $d \leq 5$. Their preprint,
which we learned of after this paper was nearly complete, gives independent proofs of the results of Section
\ref{sec:field-counting} (with different values of $\alpha_d$ and $\beta_d$).
\item 
An old theorem of Erd\H{o}s obtained the mean value of the least quadratic nonresidue to a varying modulus, or equivalently the mean value of the least non-split prime in a quadratic field. Generalizing this, Martin and Pollack \cite{MP2013} obtained the
mean value of the least prime that does not split completely in $K$, averaged over all cubic fields $K$.
\item
Zhao \cite{ZhaoHopefullySomeday} counted cubic extensions of $\mathbb{F}_q(t)$ with a power saving error term, and
the second author and Xiong \cite{TX2014} applied this to prove\footnote{Analogues of the results in Section \ref{sec:field-counting}, i.e.,
an extension of Zhao's results to allow for local conditions, have yet to be worked out in full detail in the function field setting. This is not likely to be especially difficult, but for the moment
the results of \cite{TX2014} should be considered incomplete.}
that the number of $\mathbb{F}_q$-rational
points on random trigonal covers of $\mathbb{P}^1(\mathbb{F}_q)$ converges to a Gaussian distribution, as does the number
of zeros of the associated zeta function in prescribed arcs.
(See also the bibliography of \cite{TX2014} for further related works.)
\end{itemize}
In this note we give another application and prove an Erd\H{o}s-Kac theorem for number fields. For an integer
$n$, write $\omega(n)$ for the number of prime factors of $n$ (counted without multiplicity).
The classical Erd\H{o}s-Kac theorem states that, on average over $1 \leq n \leq X$, the distribution of $\omega(n)$
approaches a Gaussian with mean and variance $\log \log X$. 

To formulate a number field analogue, write $\omega(K) := \omega(D_K)$ for a number field $K/\mathbb{Q}$; we will study the distribution
of the values $\omega(K)$. Related questions have been studied previously; for example, Belabas and Fouvry \cite{FK1999} proved
that there are infinitely many $S_3$-cubic fields whose discriminant is fundamental and has at most $7$ prime factors.\footnote{Similar results for $S_4$-quartic and $S_5$-quintic
fields, albeit with a larger bound on the number of prime factors, should follow from the results described in Section \ref{sec:field-counting}.}
Although probabilistic considerations suggest that there should be infinitely many cubic, quartic, and quintic fields of prime discriminant,
it seems that this question is beyond the reach of existing methods.

However, as we will see, an Erd\H{o}s-Kac theorem for number fields is straightforward to prove:

\begin{theorem}
\label{thm:erdos-kac}
Let $d=2,3,4,$ or $5$.  For any $z\in\mathbb{R}$, we have that
\begin{equation}\label{eq:erdos-kac}
\lim_{X\to\infty} \frac{\#\{K \in \mathcal{F}_{S_d}(X): \omega(K)-\log\log X \leq z \sqrt{\log\log X}\} }{N_{S_d}(X)} = \frac{1}{\sqrt{2\pi}} \int_{-\infty}^z e^{-t^2/2}\,dt.
\end{equation}
That is, as $K$ ranges over $S_d$ number fields of degree $d$, $\omega(K)$ becomes normally distributed.
\end{theorem}

\begin{remark}
1. The $d=2$ case requires less deep input than the other cases, and, in particular, is already in the literature.  A proof of a slightly more general result can be found, for example, in \cite{KLO2013}.

2. One might also ask for a quantitative version of Theorem \ref{thm:erdos-kac}, i.e., for a version with an explicit error term. 
We could presumably do this with more care; the results would depend on the error terms in Section \ref{sec:field-counting}
and would likely be far from optimal.
\end{remark}

Our result is stated only for $S_d$-extensions, and, for the most part, we leave open the analogous question 
for other Galois groups. However, it can be proved for $d = 3$ or $5$ that such extensions may be included
in \eqref{eq:erdos-kac} as follows.  If $G_d \neq D_4$ is a transitive permutation group on $d$ letters, $d \leq 5$, then
it is known that $N_{G_d}(X) \ll X^{1 - \delta}$ for some fixed $\delta > 0$.
Such results have been proved by a variety of authors, and all of these cases may be proved by using
the Selberg sieve (exactly following \cite{ST2014}) in conjunction with an argument along the lines of Section 3.2 of
\cite{Bhargava2010}. In the $A_5$ case, 
the (easy) details of this argument, along with a stronger error term,
will appear in forthcoming work of Bhargava, Cojocaru, and the second author \cite{BCTip}. 
In particular, for $d = 3$ or $5$, these $G_d$-fields make a negligible contribution to \eqref{eq:erdos-kac}, and it follows
that one may count all degree $d$ extensions in Theorem \ref{thm:erdos-kac} instead of just $S_d$-extensions.

In the $D_4$ case, it was proved by 
by Cohen, Diaz y Diaz, and Olivier \cite{CDO2002, Cohen2003} 
that $N_{D_4}(X) = c_{D_4} X + O(X^{3/4 + \epsilon})$. 
The proof involves a careful study of associated arithmetic Selmer groups, and yields an explicit 
representation
for the associated Dirichlet series $\sum_{K} |\Disc(K)|^{-s}$. 
However, the density $\rho_{D_4}(q)$ of field discriminants
divisible by $q$ is not obviously multiplicative in $q$, so that more care would be needed when applying the probabilistic model we use in Section \ref{sec:erdos-kac}.
We leave the question of whether an Erd\H{o}s-Kac theorem should
hold for $D_4$-extensions for later work.

For $S_3$-{\itshape sextic} extensions (i.e., extensions Galois over $\Q$ with Galois group $S_3$), 
our methods should yield a proof of the analogue of Theorem \ref{thm:erdos-kac}. The required analogue of
Theorem \ref{thm:field_count} was essentially proved in \cite{TT6}, except that the dependence of the error terms
on $q$ was not made explicit there. It would not be difficult to do so, but the details are a bit messy and so we omit them here.
\smallskip

Belabas has implemented an algorithm to quickly enumerate cubic fields (see \cite{Belabas1997}, with accompanying software available from
his website),
which allows us to compare our results for $d = 3$ to numerical data. The graph on the left gives the count of $\omega(K)$ for all cubic fields
$K$ with $|\Disc(K)| \leq 10^8$; for comparison, the right graph gives the distribution of $\omega(n)$ over $n \leq 10^8$.

\medskip
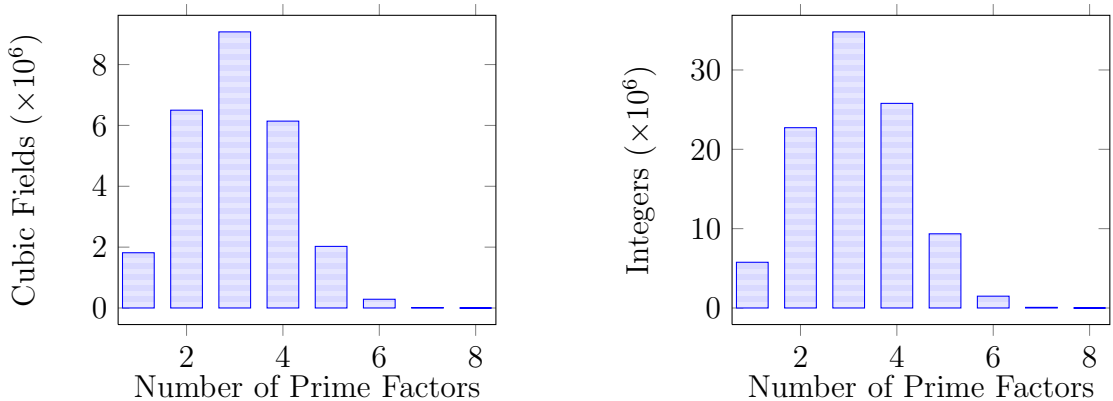
\begin{tikzpicture}
\begin{axis}[ybar,enlargelimits=0.06, 
    width=0.4\textwidth,
    legend style={at={(0.5,-0.2)},
      anchor=north,legend columns=-1},
    bar width=12pt,
xlabel={Number of Prime Factors},
ylabel={Cubic Fields ($\times 10^6$)},
]
\addplot
[draw=blue,pattern=horizontal lines light blue
] 
coordinates
{(1, 1.815920) (2, 6.501043) (3, 9.074147) (4, 6.141365) (5, 2.024511) (6, 0.287351) (7, 0.013045) (8, 0.000093)};
	
\end{axis}
\end{tikzpicture}
\hspace{0.5in}
\begin{tikzpicture}
\begin{axis}[ybar,enlargelimits=0.06, 
    width=0.4\textwidth,
    legend style={at={(0.5,-0.2)},
      anchor=north,legend columns=-1},
    bar width=12pt,
xlabel={Number of Prime Factors},
ylabel={Integers ($\times 10^6$)},
]
\addplot
[draw=blue,pattern=horizontal lines light blue
] 
coordinates
{(1, 5.762859) (2, 22.724609) (3, 34.800362) (4, 25.789580) (5,9.351293) (6,1.490458) (7, 0.080119) (8, 0.000719)};
	
\end{axis}
\end{tikzpicture}
\medskip

The graphs are extremely similar, although the cubic discriminants have slightly fewer prime divisors on average.
This does not seem too surprising, as, for example, there are no field discriminants divisible by $2$ but not $4$.

\medskip
This paper is organized as follows. In Section \ref{sec:field-counting}, we recall the necessary results from the literature on the distribution of number fields.  In fact, while essentially all the results we need are known, not all have appeared in the literature.  We therefore also sketch proofs of these results.  In Section \ref{sec:erdos-kac}, we prove Theorem \ref{thm:erdos-kac} by adapting a proof of the classical Erd\H{o}s-Kac theorem due to Billingsley.

\section*{Acknowledgments}
We would like to thank Etienne Fouvry, Ken Ono, Ari Shnidman, Arul Shankar, and Jacob Tsimerman for helpful feedback.

\section{Counting fields with power-saving error terms}
\label{sec:field-counting}

In this section we give the estimates for counting number fields needed in the proof of our main result.

Deviating slightly from the notation in the introduction, for each $d \leq 5$, write $N_d(X)$ for the counting function of degree $d$ number fields $K$ with $|D_K| < X$,
and whose Galois closure has Galois group $S_d$ or $A_d$ over $\Q$. Additionally, write $N_d(X, q)$ for the 
count of such number fields whose discriminant is divisible by $q$. The main technical result we need is the following.

\begin{theorem}\label{thm:field_count}
Let
\begin{equation}
\alpha_2 = \frac{1}{2}, \ \ \ \alpha_3 = \frac{1}{6}, \ \ \ \alpha_4 = \frac{1}{240}, \ \ \ \alpha_5 = \frac{1}{200},
\end{equation}
\begin{equation}
\beta_2 = - \frac{1}{2}, \ \ \ \beta_3 = \frac{2}{3}, \ \ \ \beta_4 = \frac{9}{10}, \ \ \ \beta_5 = 1.
\end{equation}

Then, there are constants $c_d$ and multiplicative functions $\rho_d(q)$, the latter given by
\begin{equation}\label{eq:def_rho}
\rho_d(p) :=
\left\{\begin{array}{ll}
1 - \frac{1}{1 + p^{-1}} & \text{if $d=2$}, \\
1 - \frac{1}{1 + p^{-1} + p^{-2}} & \text{if $d=3$}, \\
1 - \frac{1 }{1 + p^{-1} + 2 p^{-2} + p^{-3} } & \text{if $d=4$}, \text{ and} \\
1 - \frac{1 }{1 + p^{-1} + 2 p^{-2} + 2 p^{-3} + p^{-4} } & \text{if $d=5$}, \end{array} \right.
\end{equation}
such that for all squarefree $q$ we have
\begin{equation}
N_d(X, q) = c_d \rho_d(q) X + O(X^{1 + \epsilon - \alpha_d} q^{\beta_d}),
\end{equation}
the implied constant being absolute.
\end{theorem}
The precise values of $\alpha_d$ and $\beta_d$ given are not optimal, and for our purposes would be relevant only if we wanted
to establish a quantitative rate of convergence to the Gaussian. For improved constants, together with a very detailed exposition of how they are obtained, we refer ahead
to forthcoming work of Ellenberg, Pierce, and Wood \cite{EPW}.
(Similar results were also obtained independently by Cho and Kim \cite{CK2014_CLT}, and previously in part by Yang \cite{Yang2009}.)
Here we offer a quick proof which does not aim at optimal
values of $\alpha_d$ and $\beta_d$. 

\smallskip
When $d = 2$, the result follows by an elementary inclusion-exclusion argument.

\smallskip

\smallskip
For $d = 3$, this is a variation of the Davenport-Heilbronn theorem,
essentially first proved in work of Belabas, Bhargava, and 
Pomerance \cite{BBP2010}, and with the constants above proved in \cite{TT2013}.
The counting functions $N_3(X, q)$ also 
have a negative secondary term of order $X^{5/6}$, as proved (independently) in \cite{BST2013} and \cite{TT2013},
and with this secondary term included above we could take $\alpha_d = 2/9, \ \beta_d = 8/9$.

\smallskip
For $d = 5$ (we handle $d = 4$ last), this result was obtained for $q = 1$ 
in an elegant short paper of Shankar and Tsimerman \cite{ST2014} by means of the Selberg sieve.
To extend their result to general $q$ it is only necessary to adjust their final computation on p. 7, in which they `sieve to fields.'  We very briefly recap the situation at this point in their work, but our exposition is not to be taken independent of theirs.  Let $V_\Z$ denote the space of quadruples of integral skew-symmetric $5\times 5$ matrices, equipped with an action of $G_\Z=\mathrm{GL}_4(\Z)\times\mathrm{SL}_5(\Z)$, and recall that the $G_\Z$-orbits of $V_\Z$ parametrize quintic rings, with quintic fields corresponding to orbits of irreducible rings that are everywhere maximal.  Combining equations (2) and (8) of \cite{ST2014}, the number of orbits of discriminant up to $X$ that are \emph{i)} irreducible, but outside the `main ball,' \emph{ii)} reducible, and inside the main ball, or \emph{iii)} corresponding to an order in a non-$S_5$ quintic field, is 
\[
O(X^{199/200+\epsilon}).
\]
Thus, to obtain Theorem \ref{thm:field_count}, it remains to count those points in the main ball which are everywhere maximal and satisfy the specified local conditions at $q$.  

Let $V_{\Z, q}^{(i)}$ denote the subset of $x \in V_{\Z}^{(i)}$ corresponding to quintic rings which are maximal at $q$
and for which $q \mid \Disc(x)$, and write $k'_q$ for its density. Then, analogously to p. 7 of \cite{ST2014} and with $U_p$, $W_d$, and $k_d$ as they are there, we find that the desired number of points in the main ball is
\begin{align*}
N_{12}^*(\cap_p U_p \cap V_{\Z, q}^{(i)}, X) = &
\sum_{d \in \mathbb{N}, \ (d, q) = 1} \mu(d) N_{12}^*(W_d \cap V_{\Z, q}^{(i)}, X)
\\
= & 
\sum_{d < T, \ (d, q) = 1} \Big( c_i \mu(d) k_d k'_q X + O(X^{39/40} d^{\epsilon} q^{1 + \epsilon}) \Big)
+ \sum_{d > T} O_{\epsilon} (X/d^{2 - \epsilon})
\\
= & 
\sum_{d \in \mathbb{N}, \ (d, q) = 1} c_i \mu(d) k_d k'_q X + O_{\epsilon}(X/T^{1 - \epsilon} + X^{39/40} T^{1 + \epsilon} q^{1 + \epsilon})
\\
= & \
c_i k'_q \prod_{p \nmid q} (1 - k_p) X
+ O_{\epsilon}(X/T^{1 - \epsilon} + X^{39/40} T^{1 + \epsilon} q^{1 + \epsilon}).
\end{align*}
In the second line we have used that  
$W_d \cap V_{\Z, q}$
consists of the union of  
$O(d^{78 + \epsilon} q^{79 + \epsilon})$ translates of
$(dq)^2 V_\Z$, with $(Tq)^2 < X^{1/40}$ so that we may 
apply (4) of \cite{ST2014}, and
choosing $T = X^{1/80} q^{-1}$ we obtain an error term of $X^{79/80 + \epsilon} q^{1 + \epsilon}$.  Combining this in the most simple-minded way with the previous error $O(X^{199/200+\epsilon})$, we obtain Theorem \ref{thm:field_count}.


\begin{remark} This corrects a typo in \cite{ST2014}, where $N(-)$ should be replaced with $N_{12}^*(-)$ on p. 7, as we did above.
(In other words, we count elements of $V_{\Z}$ without regard to reducibility but with the requirement that $a_{12} \neq 0$; then (2) and (8)
of \cite{ST2014} establish the stated bounds.)
\end{remark}

\medskip
For $d = 4$, it was proved by Bhargava \cite{Bhargava2004} that $G_\Z := \GL_3(\Z) \times \GL_2(\Z)$-orbits 
on $V_\Z = (\Sym^2 \Z^3 \otimes \Z^2)^*$ are in bijection with isomorphism classes of pairs $(Q, R)$, where $Q$ is a quartic ring and
$R$ is a {\itshape cubic resolvent} of $Q$. Bhargava used this parameterization \cite{Bhargava2005} to obtain the asymptotic
density of quartic fields, and with Belabas and Pomerance \cite{BBP2010} refined this with a power-saving error term.

An analogue of (4) of \cite{ST2014} is given
in Theorem 4.11 of \cite{BBP2010}. For any $G_{\Z}$-invariant subset $S \subseteq V_{\Z}$, define
$N_{11}^*(S, X)$ to be the average number of elements $v \in S$ in a fundamental domain
for $G_{\Z} \backslash V$, with $a_{11} \neq 0$ and discriminant less than $X$, as in Theorem 4.11 of
\cite{BBP2010} or Sections 2 and 3 of \cite{ST2014}.  Theorem 4.11 then
yields that for any translate $L$
of $m V_{\Z}$, for $m \leq X^{1/12}$, we have
\begin{align}\label{eqn:bbp_4}
N_{11}^*(L \cap V_{\Z}^{(i)}, X) &
= c_i m^{-12} X + O\big( m^{-12} X^{23/24} + m^{-7} X^{11/12} \big)
\\
&= c_i m^{-12} X + O\big( m^{-7} X^{23/24}\big).
\end{align}

Imitating the $d = 5$ argument given above (alternatively, we could follow \cite{BBP2010}) we find the desired number of points in the main ball is
\begin{align*}
N_{11}^*(\cap_p U_p \cap V_{\Z, q}^{(i)}, X) = &
\sum_{d \in \mathbb{N}, \ (d, q) = 1} \mu(d) N(W_d \cap V_{\Z, q}^{(i)}, X)
\\
= & 
\sum_{d < T, \ (d, q) = 1} \Big( c_i \mu(d) k_d k'_q X + O(X^{23/24} d^{8 + \epsilon} q^{9 + \epsilon}) \Big)
+ \sum_{d > T} O_{\epsilon} (X/d^{2 - \epsilon})
\\
= & 
\sum_{d \in \mathbb{N}, \ (d, q) = 1} c_i \mu(d) k_d k'_q X + O_{\epsilon}(X/T^{1 - \epsilon} + X^{23/24} T^{9 + \epsilon} q^{9 + \epsilon})
\\
= & \
c_i k'_q \prod_{p \nmid q} (1 - k_p) X
+ O_{\epsilon}(X/T^{1 - \epsilon} + X^{23/24} T^{9 + \epsilon} q^{9 + \epsilon}).
\end{align*}
Similarly to above, 
we have used that  
$W_d \cap V_{\Z, q}$
consists of the union of  
$O(d^{22 + \epsilon} q^{23 + \epsilon})$ translates of
$(dq)^2 V_\Z$; the bound quoted for $d > T$ follows from Lemma 4.3 of \cite{BBP2010}.
We then choose $T = X^{1/240} q^{-9/10}$, observe that the condition $(Tq)^2 < X^{1/12}$ required for \eqref{eqn:bbp_4} is easily
satisfied, and obtain an error term $\ll X^{239/240} q^{9/10}$. The main term is as in \cite{Bhargava2005}. 

Finally,
Lemmas 4.9 and 4.10 of \cite{BBP2010} establish that up to an error of $O(X^{11/12 + \epsilon})$ 
we may exchange the `main ball' condition on $a_{11}$ for `total irreducibility', which, by Theorem 4.1 of \cite{BBP2010}), is equivalent to the property that the associated quartic ring is an $S_4$-quartic order. This completes the proof.

\begin{remark} The tail estimate used to estimate the sum over $d > T$ (Lemma 4.3 of \cite{BBP2010}) is proved for $S_4$-quartic orders rather
than for elements of $V_\Z$, and is therefore subject to the $O(X^{11/12 + \epsilon})$ error term from Lemmas 4.9 and 4.10 of
\cite{BBP2010} described above. Since any given order or element of $V_\Z$ is counted at most $O(X^{\epsilon})$ times (at most once for each $d > T$
dividing the discriminant), this error term may be absorbed into the stated error term.
\end{remark}

\section{The Erd\H{o}s-Kac machinery}
\label{sec:erdos-kac}

In order to prove Theorem \ref{thm:erdos-kac}, we adapt a proof of the classical Erd\H{o}s-Kac theorem due to Billingsley \cite{Billingsley1974}.  We proceed via the method of moments, which in this case is made easier as we are content with establishing a qualitative result; in particular, the moment calculation need not be uniform.  With not too much more effort, our results could be made uniform, and therefore quantitative, but in order to obtain an optimal result, we would likely need to proceed in a different direction.

Fix $d=2,3,4,$ or $5$, and, for each squarefree $q$, define $\rho_d(q)$ as in Theorem \ref{thm:field_count}.
Define a random variable $R_{d,p}$ to be $1$ with probability $\rho_d(p)$ and 0 with probability $1-\rho_d(p)$, so that
$R_{d,p}$ models the event that $p$ ramifies in a number field of degree $d$.  As a consequence of the central limit theorem (e.g., via Lyapunov's criterion), we know that the quantity
\[
R_d(Z) := \sum_{p\leq Z} R_{d,p}
\]
becomes normally distributed as $Z\to\infty$ and, provided $Z=X^\delta$ for some $\delta>0$, with mean and variance each $\mu(X):=\log\log X$.  Thus, for fixed $k$, we have that
\begin{equation}\label{eq:rd_normal}
\mathbb{E}\left( (R_d(Z)-\mu(X))^k \right) = c_k \mu(X)^{k/2} +o_k(\mu(X)^{k/2}),
\end{equation}
where
\[
c_k = \left\{ \begin{array}{ll} \frac{k!}{2^{k/2}(k/2)!}, & \text{if $k$ is even, and} \\ 0, & \text{if $k$ is odd} \end{array} \right.
\]
is the $k$-th moment of the standard Gaussian.

Following Billingsley's proof and the notation of Section \ref{sec:field-counting}, the key idea is to compare the $k$-th moment
\[
M_{d,k}(X) := \frac{1}{N_d(X)} \sum_{K\in\mathcal{F}_d(X)} \left(\omega(K)-\mu(X)\right)^k
\] 
of the quantity we are interested in to that of the model; here, $\mathcal{F}_d(X):=\{K/\mathbb{Q}: [K:\mathbb{Q}]=d, \mathrm{Gal}(K/\mathbb{Q})\cong S_d \text{ or } A_d, \text{ and } |D_K|\leq X\}$, so that $N_d(X)=\#\mathcal{F}_d(X)$.  In particular, as the moments of the Gaussian determine the distribution (for example, the moment generating function converges on all of $\mathbb{C}$, and the distribution function can be recovered simply by Fourier inversion), our main theorem follows if we prove the analogue of
\eqref{eq:rd_normal} for 
each of these moments.

We let $Z=X^{\alpha_d / 2 k (\beta_d + 1)}$,
write
$\omega(K;Z):=\#\{p\mid D_K : p \leq Z\}$, and for each number field $K$ of degree $d$ and each prime $p$, define 
\[
\delta_p(K) := \left\{\begin{array}{ll} 1, & \text{if } p\mid D_K, \text{ and } \\ 0, & \text{otherwise.} \end{array} \right.
\]
Finally, for any $k^\prime\leq k$, define
\[
\tilde{M}_{d,k^\prime}(X,Z) := \frac{1}{N_d(X)} \sum_{K\in\mathcal{F}_d(X)} \omega(K;Z)^{k^\prime}.
\]
We now observe that
\begin{eqnarray*}
\tilde{M}_{d,k^\prime}(X,Z)
	&=& \frac{1}{N_d(X)} \sum_{K\in\mathcal{F}_d(X)} \left(\sum_{p\leq Z} \delta_p(K) \right)^{k^\prime} \\
	&=& \sum_{p_1,\dots,p_{k^\prime} < Z} \frac{1}{N_d(X)} \sum_{\begin{subarray}{c} K\in\mathcal{F}_d(X): \\ p_1,\dots,p_k^\prime\mid D_K \end{subarray}} 1 \\
	&=& \sum_{p_1,\dots,p_{k^\prime} < Z} \left(\rho_d(\mathrm{lcm}(p_1,\dots,p_{k^\prime})) + O\left(X^{-\alpha_d}Z^{k^\prime \beta_d}\right)\right) \\
	&=& \mathbb{E}(R_d(Z)^{k^\prime}) + O\left(X^{-\alpha_d/2}\right),
\end{eqnarray*}
by Theorem \ref{thm:field_count}, the construction of the random variables $R_{d,p}$, and the choice of $Z$.  Letting $M_{d,k}(X,Z):=\mathbb{E}((\omega(K;Z)-\mu(X))^k)$, this shows that 
\begin{eqnarray*}
M_{d,k}(X,Z)
	&:=& \frac{1}{N_d(X)} \sum_{K\in\mathcal{F}_d(X)} \left(\omega(K;Z)-\mu(X)\right)^k \\
	&=& \sum_{j=0}^k \left(k \atop j\right) (-\mu(X))^j \tilde{M}_{d,k-j}(X,Z) \\
	&=& \sum_{j=0}^k \left(k \atop j\right) (-\mu(X))^j \left(\mathbb{E}(R_d(Z)^{k-j}) + O(X^{-\alpha_d/2})\right) \\
	&=& \mathbb{E}((R_d(Z)-\mu(X))^k) + O_k(X^{-\alpha_d/2} \mu(X)^{k}).
\end{eqnarray*}
In particular, by \eqref{eq:rd_normal}, we see that $M_{d,k}(X,Z) = c_k \mu(X)^{k/2} +o_k(\mu(X)^{k/2})$, where we recall that $c_k$ is the $k$-th moment of the standard Gaussian.

Finally, it remains to compare the truncation $M_{d,k}(X,Z)$ to the full moment $M_{d,k}(X)$.  For this, we see that
\begin{eqnarray*}
M_{d, k}(X) &=& \frac{1}{N_d(X)} \sum_{K\in\mathcal{F}_d(X)} \left(\omega(K)-\mu(X)\right)^k \\
	&=& \frac{1}{N_d(X)} \sum_{K\in\mathcal{F}_d(X)} \left(\omega(K;Z)-\mu(X) + O(k)\right)^k \\
	&=& M_{d,k}(X,Z) + O_k\left(\sum_{j=0}^{k-1} \frac{1}{N_d(X)} \sum_{K\in\mathcal{F}_d(X)} |\omega(K;Z)-\mu(X)|^j\right) \\
	&=& M_{d,k}(X,Z) + o_k\left(\mu(X)^{k/2}\right),
\end{eqnarray*}
where in advancing to the final line we employ the Cauchy-Schwarz inequality for those $j$ that are odd.  This completes the proof.

\bibliographystyle{alpha}
\bibliography{discriminant_ek}

\end{document}